\documentclass[leqno,12pt]{article} 
\newcommand {\s}[1]{{\cal #1}} \newcounter{sectiune}
\newcounter{lemma} 
\newenvironment{sectiune}{\setcounter{lemma}{0}%
                         \setcounter{equation}{0}%
                         \stepcounter{sectiune}%
                         \bigskip \noindent\large \bf%
                         \arabic{sectiune}.  }%
                         \bigskip

\newenvironment{lemma}{\bigskip
                      \stepcounter{lemma}\begingroup \noindent
                        \bf Lemma
                       \arabic{sectiune}.\arabic{lemma}. \it}{
                        \endgroup\vskip\baselineskip}

\newenvironment{remark}{\bigskip
                      \stepcounter{lemma}\begingroup \noindent
                      \it Remark
                       \arabic{sectiune}.\arabic{lemma}. \rm}{
                        \endgroup\vskip\baselineskip}

\begin{document}


\noindent
{\bf \Large Holomorphic  vector bundles on\\ primary Kodaira
surfaces}

\bigskip

\noindent
{Marian Aprodu, Vasile Br\^{\i}nz\v{a}nescu and Matei Toma}

\bigskip

\noindent
{\footnotesize
Institute of Mathematics of the Romanian Academy\\
P.O.Box 1-764, RO-70700, Bucharest, Romania\\
E-mail: Marian.Aprodu@imar.ro\\
E-mail: Vasile.Brinzanescu@imar.ro\\
E-mail: Matei.Toma@imar.ro\\
Universit\"at Osnabr\"uck, Fachbereich Mathematik/Informatik,\\
D-49069, Osnabr\"uck, Germany \\
E-mail: matei@mathematik.Uni-Osnabrueck.DE}

\bigskip

\bigskip

\noindent{\large \bf
Introduction}

\bigskip

\noindent
Let $X$ be a smooth compact complex surface. A classical problem
is to decide which topological complex 
vector bundles admit holomorphic structures, or equivalently,
to find all the triples 
$(r,c_1,c_2)\in {\bf N}\times NS(X)\times {\bf Z}$,
$r\geq 2$ for which there exists a rank $r$ holomorphic
vector bundle $\cal E$ on $X$ with Chern classes
$c_1({\cal E})=c_1$ and $c_2({\cal E})=c_2$.

For projective surfaces, Schwarzenberger (cf. \cite{Schw})
proved that any triple $(r,c_1,c_2)\in {\bf N}\times NS(X)\times {\bf
Z}$,
$r\geq 2$ comes from a rank $r$ holomorphic vector bundle.
In contrast to this situation, for
non-projective surfaces, there is a natural
{\em necessary} condition
for the existence problem (cf. \cite{BaL} Theorem 3.1
for the general case;
for the rank-2 case see \cite{BrF} Proposition 1.1, \cite{ElFo}):

$$
\Delta (r,c_1,c_2):=\frac{1}{r}\left( c_2-\frac{r-1}{2r}c_1^2
\right) \geq 0.
$$

Using essentially extensions of coherent sheaves 
one proves the following
{\em sufficient}
condition 
for the existence of holomorphic vector bundles
(cf \cite{BaL}; see also \cite{BrF}, \cite{LeP}):
$$
\Delta (r,c_1,c_2)\geq m(r,c_1),
$$
where
$$
m(r,c_1):=-\frac{1}{2r}max\left\{ 
\sum\limits_{i=1}^r\left( \frac{c_1}{r}-\mu _i\right) ^2,
\; \mu _1,...,\mu _r\in NS(X),\sum\limits_{i=1}^r\mu _i
=c_1\right\} ,
$$
with the sole excepted case: $X$ a $K3$ surface
with algebraic dimension zero, $c_1$ divisible
by $r$ in $NS(X)$ and $\Delta \big( r,c_1,c_2\big) =\frac{1}{r}$.
This  result  also has a converse, 
which
is (see \cite{BaL}, \cite{BrF}, \cite{LeP}) any {\em filtrable} rank $r$ 
holomorphic vector bundle $\cal E$ on a non-algebraic surface $X$ 
with Chern classes $c_1({\cal E})=c_1$ and $c_2({\cal E})=c_2$
satisfies the inequality
$$\Delta ({\cal E}):=\Delta (r,c_1,c_2)\geq m(r,c_1).$$

Therefore the only unknown situations are
in the range
$\Delta (r,c_1,c_2)\in \left[ 0, m(r,c_1)\right) $.
If $m(r,c_1)\neq 0$, 
this interval is non-empty,
and 
in order to solve the existence problem, one
has to construct holomorphic vector bundles
having the corresponding discriminant $\Delta $
lower that $m(r,c_1)$. Of course, all these
vector bundles will be {\em non-filtrable} 
and the difficulty of the problem resides in the lack of a general 
construction method in this case
(for more details see, for 
example, \cite{Bri3}).
One is compelled therefore to 
focus on
particular classes
of surfaces and find some specific construction methods.

The main result in the present paper is the following:\\\\
{\bf Theorem.}
{\em Let $X$ be a primary Kodaira
surface. 
Then for
a triple $(r,c_1,c_2)\in {\bf N}\times NS(X)\times {\bf Z}$,
with $r\geq 2$,
there exists a holomorphic rank-$r$ vector bundle
${\cal E}$ on $X$ such that 
$c_1({\cal E})=c_1$ and
$c_2({\cal E})=c_2$
if and only if $\Delta (r,c_1,c_2)\geq 0$.}

\bigskip

From now on
by a {\em vector bundle} we shall mean a
{\em holomorphic vector bundle}, and
a {\em curve} will always be a smooth,
complex projective curve.

Recall that a primary Kodaira surface $X$ is a principal elliptic bundle
over some elliptic curve $B$.
The case 
when $r \Delta(r,c_1,c_2)$ is an integer 
has been
previously solved (cf. \cite{To2}; see also \cite{To1}) by using
{\em unramified} coverings of $B$ 
by
another suitable elliptic curve and suitable deformations
of sheaves.
For 
the remaining
cases,  we use 
{\em curves of genus two with elliptic 
differentials}, i.e. curves $C$ of genus 2 which admit a non-constant 
morphism to an elliptic curve, which 
does not factor over an isogeny of 
the
elliptic curve. 
These curves were studied 
extensively by Bolza, Humbert, Picard, Poincar\'{e} and many others 
(see Chapter XI of Krazer's book \cite{Kr}), and a
remarkable thing about them is that the covers
occur in pairs. More precisely (for details, see for
example \cite{FrKa}), if we
consider $f:C\rightarrow B$ such a covering, and
choose $F=\mbox{Ker}(J_C\rightarrow B)$,
then there is a {\em complementary} natural
covering $g:C\rightarrow F$
(induced by the embedding of $C$ in $J_C$), giving
rise to an isogeny of degree $r^2$,
$f^*\times g^*:B\times F\rightarrow J_C$. 
Its kernel $H$ is the graph of an isomorphism
$\psi :B[r]\rightarrow F[r]$ (which 
is {\em anti-isometric} with respect to the Weil pairings).

Conversely, one can 
get such curves out of a particular case
of the "Basic Construction" and "Reducibility Criterion"
(cf. \cite{Ka2} and \cite{We}; see also \cite{Ka1},
\cite{FrKa}). Suppose
one starts with $B$ and $F$  two elliptic curves, and
$\Theta =0\times F+B\times 0$ the canonical principal polarization
on $B\times F$ given by the product structure. Assume  
$\psi :B[r]\rightarrow F[r]$ is
an isomorphism between the $r$-torsion subgroups,
anti-isometric with respect to the Weil pairings.
Then, by denoting $H_\psi =\mbox{Graph}(\psi )$
(which is a maximal isotropic subgroup
of $(B\times F)[r]$),
and $p:B\times F\rightarrow (B\times F)/H_\psi $
the canonical projection, it turns out that
the surface $(B\times F)/H_\psi $
carries a principal polarization
$\Theta _\psi $ such that $p^*\Theta _\psi $
is linear equivalent to $r\Theta $
(cf. \cite{FrKa}, \cite{Ka2}; see also \cite{Mumf1}).
If the anti-isometry $\psi$ is irreducible
(cf. \cite{We}; see  \cite{Ka2} for an
analysis of this case), then  
the linear
system $|\Theta _\psi |$ contains a smooth curve $C$
of genus 2.
In particular, the principally
polarized abelian surface $\big( (B\times F)/H_\psi ,\Theta _\psi \big)$
is isomorphic to the Jacobian of $C$, and 
$C$ covers $r$-to-one
the curves $B$ and $F$.

\bigskip

{\em Acknowledgements}. 
The first named author was supported
by a DFG post-doctoral fellowship within the
Graduate Programme "Complex Manifolds" at the University
of Bayreuth, and by a visiting fellowship at the
Abdus Salam International Centre for 
Theoretical Physics (ICTP) in Trieste.
The author expresses his special thanks to the above
mentioned institutions for hospitality during the 
preparation of this paper.

This work has been  done while the 
second named author was visiting Kaiserslautern University with a 
resumption of the Humboldt stipend, and
ICTP Trieste as a Regular Associate.  He is grateful to G. Trautmann, to 
the Alexander von Humboldt-Stiftung, and to
ICTP for this support.

\begin{sectiune}
{\bf Topological information on principal elliptic bundles}
\end{sectiune}

\noindent
This section is devoted to some results
which are used in the proof of the Theorem.
We start by giving a formula for computing 
the degree of a covering of an elliptic curve.

\begin{lemma}
Let $\gamma :B\rightarrow E$ be a covering of an elliptic
curve with a curve of genus $g$. Choose
$\{ \alpha _1,...,\alpha _g,\beta _1,...,\beta _g\} 
\in H_1(B,{\bf Z})$ a symplectic basis.
If $\gamma _*:H_1(B,{\bf Z})\rightarrow H_1(E,{\bf Z})$
denotes the natural push-forward map, then
$$\mbox{\em deg
}(\gamma )=\sum\limits_{i=1}^g
\gamma _*(\alpha _i). \gamma _*(\beta _i).$$
\end{lemma}

\noindent
{\em Proof.} Consider the Abel-Jacobi embedding
of $B$ in its Jacobian $J_B$, corresponding to a
chosen point in $B$. Following classical notations
(cf. \cite{LanBi}, p. 328), the class of $B$ as a cycle 
in $J_B$ is represented by $\{ \widetilde{W_1}\}
\in H_2(J_B,{\bf Z})$, $\{ \widetilde{W_1}\}=
-\sum\limits_{i=1}^g\alpha _i\star \beta _i$,
where "$\star $" denotes the Pontrjagin product.
If we denote 
$V=H^0(B,K_B)^*$, then 
for any divisor $D$ on $J_B$ whose Chern class
is represented, by means of Appell-Humbert theorem,
by an alternating form $A:V\times V\rightarrow {\bf R}$,
we have (cf. \cite{LanBi}, p.42-43, p. 104, p.106):
\begin{equation}
\left( \{ \widetilde{W_1}\} .\{ D\} \right) =\sum\limits_{i=1}^g
 A(\alpha_i,\beta_i). 
\end{equation}

If one considers now $D$ a fiber 
of the norm map $Nm(\gamma ):J_B\rightarrow E$,
then the associated alternating form is
$A:V\times V\rightarrow {\bf R}$, defined
by $A(u,v)=\gamma _*(u).\gamma_*(v)$ for
all $u,v\in H_1(B,{\bf Z})$.
Now apply the relation (1.1), and 
observe that $\mbox{deg}(\gamma )
=\left( \{ \widetilde{W_1}\} .\{ D\} \right)$,
and the proof of the Lemma 
is over.

\begin{remark}
In particular, if we also choose $\{ \varepsilon _1,\varepsilon _2\} 
\in H_1(E,{\bf Z})$ a symplectic basis, and the matrix
of $\gamma _*$ in the two given bases is
$$\gamma_*=\left( \begin{array}{cccccc}
a_{11} & ... & a_{1g} & b_{11} & ... & b_{1g} \\
a_{21} & ... & a_{2g} & b_{21} & ... & b_{2g}
\end{array}\right) ,$$
then it easily follows that
\begin{equation}
\mbox{deg}(\gamma )=\sum\limits_i \left|
\begin{array}{cc}
a_{1i} & b_{1i} \\
a_{2i} & b_{2i}
\end{array}
\right|
.
\end{equation}
\end{remark}

Consider next a curve $B$  of genus $g\geq 1$, $E$ an elliptic curve, 
and $X\stackrel{\pi}{\longrightarrow}B$
a principal elliptic bundle over $B$ with
fiber $E$ (cf. \cite{Kod}, \cite{Hf}, \cite{Bri3} 
for precise definitions).
If we suppose moreover, that 
$X\stackrel{\pi }{\longrightarrow}B$ is not
topologically trivial, 
then we have the following isomorphisms
compatible with the cup-products
(cf. \cite{Bri1}, \cite{Bri2}):
\begin{equation}
H^2(X,{\bf Z})/\mbox{Tors}\big( H^2(X,{\bf Z})\big)
\cong H^1(B,{\bf Z})\otimes _{\bf Z}H^1(E,{\bf Z}),
\end{equation} 
and
\begin{equation}
NS(X)/\mbox{Tors}\big( NS(X)\big)
\cong \mbox{Hom}(J_B,E^\vee),
\end{equation}
where $J_B$ denotes the Jacobian variety of $B$ and
$E^\vee$ is the dual curve of $E$.
The torsion of $H^2(X,{\bf Z})$ (as well as of $NS(X)$)
is generated by the class of a fiber. 

To avoid confusions, throughout the rest of the paper, for
an element $c\in H^2(X,{\bf Z})$, we will denote
by $\widehat{c}$ its class modulo $\mbox{Tors}\big( H^2(X,{\bf Z})
\big)$.

We prove next a formula for
computing the intersection form on $H^2(X,{\bf Z})$.
It suffices, of course, to make it explicit for
classes in $H^2(X,{\bf Z})/\mbox{Tors}\big( H^2(X,{\bf Z})\big)$.
Consider then $\{ \alpha_1,...,
\alpha_g,\beta_1,...,\beta_g\} \subset H_1(B,{\bf Z})$
a symplectic basis with Kronecker duals
$\alpha_1^*,...,\alpha_g^*,
\beta_1^*,...,\beta_g^*\in H^1(B,{\bf Z})$,
and $\{ e_1,e_2\}\subset H^1(E,{\bf Z})$ a symplectic 
basis. 
Then  any class
$\widehat{c} \in H^2(X,{\bf Z})/\mbox{Tors}\big( H^2(X,{\bf Z})\big)$
can be expressed as
$$
\widehat{c}=\sum\limits_{k=1}^2\sum\limits_{i=1}^g 
\big( a_{ki}(\alpha_i^*\otimes e_k)
+b_{ki}(\beta_i^*\otimes e_k)\big) ,
$$
and the self-intersection computes as
$$
\widehat{c}^2=-2\sum\limits_i \left( a_{1i}b_{2i}-b_{1i}a_{2i}\right) .
$$

Let now $\widehat{c}\in NS(X)/\mbox{Tors}\big( NS(X)\big)$. 
Choosing a base-point in $B$,
we can think of  the cohomology
class $\widehat{c}$ as being  a covering map
$\widehat{c}:B\rightarrow E$.
By the considerations above, we
get the following result
(for a proof using the Riemann-Roch Theorem
see \cite{Te}, Remark 1.11):

\begin{lemma}
With the previous notations, 
$\widehat{c}^2=-2\mbox{\em deg}\left(  \widehat{c}\right) .$
\end{lemma}

Suppose now that $C\stackrel{f}{\longrightarrow}B$
is a covering;
in this case, $Y=X\times _BC\rightarrow C$ is
a principal elliptic bundle 
and $Y\rightarrow C$ is topologically
non-trivial as soon as $X\rightarrow B$ is topologically
non-trivial.

We give next a description
of the push-forward map induced between the cohomology
groups by the natural covering map $Y\stackrel{\varphi }
{\longrightarrow }X$. Firstly, let us remark that  the class of
a fiber in $Y$ pushes-forward to the class of
a fiber in $X$, 
therefore, we only need to describe the
push-forward map $\varphi_*$ induced 
between $H^2(Y,{\bf Z})/\mbox{Tors}\big( H^2(Y,{\bf Z})\big)$
and $H^2(X,{\bf Z})/\mbox{Tors}\big( H^2(X,{\bf Z})\big)$.
Seeing
the elements in 
$H^2(Y,{\bf Z})/\mbox{Tors}\big( H^2(Y,{\bf Z})\big)$
as maps from $H_1(C,{\bf Z})$ to $H_1(E^\vee ,{\bf Z})\cong
H^1(E,{\bf Z})$,  and
identifying the pull-back morphism 
$f^*:J_B\rightarrow J_C$ with its rational representation,
we have the following:

\begin{lemma}
With the previous notations, $\varphi _*(\widehat{c})=
\widehat{c}\circ f^*$, 
for any class
$\widehat{c}\in H^2(Y,{\bf Z})/\mbox{\em Tors}\big(
H^2(Y,{\bf Z})\big)$.
\end{lemma}

\noindent
{\em Proof.} The morphism $f^*:H_1(B,{\bf Z})
\rightarrow H_1(C,{\bf Z})$ is the pull-back
morphism in homology via Poincar\' e duality.
By using the canonical isomorphism $(1.3)$,
we get for any element $z\in H_1(B,{\bf Z})$ the following
equalities:
$$
\varphi_*(\widehat{c})(z)=\big( (f_*\otimes id)(\widehat{c})\big) /z,
\mbox{ and } (\widehat{c}\circ f^*)(z)=\widehat{c}/f^*(z),
$$
where $"/"$ is the {\em slant product}.
 Since $f_*$ and $f^*$ are transposed to each other,
we conclude.

\begin{sectiune}
Vector bundles on primary Kodaira surfaces via coverings
\end{sectiune}

\noindent
Suppose that $r\geq 2$ is an integer,
and $X\stackrel{\pi}{\rightarrow} B$ is
a primary Kodaira surface over the elliptic 
curve $B$, with fiber $E$. Let $C$ be a curve of genus
$2$, and  let $C\stackrel{f}{\rightarrow} B$ be
a (ramified) covering of degree $r$. Set 
$Y=X\times _BC\rightarrow C$, which is a
principal elliptic bundle over $C$, and covers $r$-to-$1$
$X$ by the natural map $Y\stackrel{\varphi}{\rightarrow} X$.

If one considers now a line bundle
$L$ on $Y$, the push-forward sheaf
${\s E}=\varphi _*L$ is actually a rank-$r$ vector bundle
on $X$ with
\begin{eqnarray}
&&c_1({\cal E})\equiv \varphi_*c_1(L)\mbox{ mod Tors}\big( NS(X)\big) ,
\\
\nonumber
&&
\Delta ({\cal E})=\frac{1}{2r^2}
\left[ (\varphi_*c_1(L))^2-rc_1^2(L)\right] .
\end{eqnarray}
(The push forward $\varphi_*$ on cohomology is obtained as usual via Poincar\'e
duality on both source and target space.)

Moreover, 
one can easily see that for any two line bundles
$L$ and $L'$ on $Y$,
the following holds:
\begin{equation}
c_1(\varphi_*(L\otimes L'))=c_1(\varphi_*L)+\varphi_*c_1(L').
\end{equation}

If the covering $f$ does not factor through an isogeny of $B$
we can give {\em algebraic} interpretations for 
$c_1({\cal E})$ and $\Delta({\cal E})$ as follows.
In this case,
$F:=\mbox{Ker}(J_C\rightarrow B)$ is an elliptic
curve;
consider $g:C\rightarrow F$ the complementary
covering (see, for example \cite{Ka2}, \cite{FrKa},
\cite{Mumf2}). Then the pull-back maps $f^*:
B\rightarrow J_C$ and $g^*:F\rightarrow J_C$
turn out to be injective, and they
give rise 
to
an isogeny $f^*\times g^*:B\times F\rightarrow J_C$.
Therefore, we can 
write $J_C\cong (B\times F)/H$, where 
$H=\mbox{Ker}(f^*\times g^*)$; 
the subgroup $H\subset (B\times F)[r]$ is in fact 
the graph of a  isomorphism from
$B[r]$ to $F[r]$, which is anti-isometric with respect
to the Weil pairings.

\begin{lemma}
In    the hypotheses above,   for  any line bundle $L$
on $Y$, ${\cal E}=\varphi _*L$ is a rank-$r$
vector bundle on $X$ with 
$\widehat{c_1({\cal E})} =  \widehat{c_1(L)}\circ 
f^*\; \mbox{ and }
r^2\Delta ({\cal E})=
\mbox{\em deg}(\widehat{c_1(L)}\circ g^*).$
\end{lemma}

\noindent
{\em Proof.} 
Let $\widehat{c}\in NS(Y)/\mbox{Tors}\big( NS(Y)\big)$ 
and $\widehat{c}_1\in NS(X)/
\mbox{Tors}\big( NS(X)\big)$ the 
classes  associated to 
$c_1(L)\in NS(Y)$ and $c_1({\cal E})\in NS(X)$. Then 
Lemma 1.4 and the formulae
(2.1) read
$\widehat{c}_1=\widehat{c}\circ f^*$ and
\begin{equation}
2r^2\Delta ({\cal E})=
\widehat{c}_1^2-r\widehat{c}^2  .
\end{equation}

We have the following diagram:
\begin{center}
\setlength{\unitlength}{1mm}
\begin{picture}(100,40)
\put(27,15){\makebox(10,10){$F$}}
\put(36,19){\vector(1,0){10}}
\put(46,21){\vector(-1,0){10}}
\put(45,15){\makebox(10,10){$J_C$}}
\put(54,21){\vector(1,0){10}}
\put(64,19){\vector(-1,0){10}}
\put(63,15){\makebox(10,10){$B$}}
\put(45,32){\vector(-1,-1){8}}
\put(45,31){\makebox(10,10){$C$}}
\put(50,32){\vector(0,-1){8}}
\put(50,16){\vector(0,-1){8}}
\put(55,32){\vector(1,-1){8}}
\put(36,16){\vector(1,-1){8}}
\put(64,16){\vector(-1,-1){8}}
\put(45,-1){\makebox(10,10){$E$}}
\put(39,30){$g$}
\put(60,30){$f$}
\put(46,27){$j_C$}
\put(47,12){$\widehat{c}$}
\put(42,23){$g_*$}
\put(55,23){$f_*$}
\put(42,16){$g^*$}
\put(55,16){$f^*$}
\end{picture}
\end{center}
where 
$j_C$ is the embedding of $C$ into its
Jacobian $J_C$, $g=g_*\circ j_C$, $f=f_*\circ j_C$,
$g_*\circ g^*=r\mbox{id}_F$,
$f_*\circ f^*=r\mbox{id}_B$, $f_*\circ g^*=0$ and
$g_*\circ f^*=0$.

%
Using Lemma 1.3 and 
(2.3)
we 
see that the formula $r^2\Delta ({\cal E})=
\mbox{deg}(\widehat{c_1(L)}\circ g^*)$ is equivalent to:
\begin{equation}
r\mbox{deg}(\widehat{c}\circ j_C)=\mbox{deg}(\widehat{c}\circ 
f^*)+\mbox{deg}(\widehat{c}\circ g^*).
\end{equation}

Let $\{W_1\}$, $\{f^*B\}$, $\{g^*F\}\in$$H^2(J_C,{\bf Z})$ be the
classes of the divisors $j_C(C)$, $f^*B$ and $g^*F$ respectively,
on the Jacobian $J_C$ of $C$. Denote by $S$ a fiber of the morphism
$\widehat{c}:J_C\rightarrow E$. Then
%
we have:
$$
\mbox{deg}(\widehat{c}\circ j_C)=\big( \{W_1\} . \{S\} \big) ,\;
\mbox{deg}(\widehat{c}\circ f^*)=\big( \{f^*B\} . \{S\} \big) ,\;
\mbox{deg}(\widehat{c}\circ g^*)=\big( \{g^*F\} . \{S\} \big) .
$$

Now, consider the diagram:
\begin{center}
\setlength{\unitlength}{1mm}
\begin{picture}(100,40)
\put(27,15){\makebox(10,10){$F$}}
\put(36,20){\vector(1,0){7}}
\put(45,15){\makebox(10,10){$F\times B$}}
\put(64,20){\vector(-1,0){7}}
\put(63,15){\makebox(10,10){$B$}}
\put(45,32){\vector(-1,-1){8}}
\put(45,31){\makebox(10,10){$C$}}
\put(50,32){\vector(0,-1){8}}
\put(50,16){\vector(0,-1){8}}
\put(55,32){\vector(1,-1){8}}
\put(36,16){\vector(1,-1){8}}
\put(64,16){\vector(-1,-1){8}}
\put(45,-1){\makebox(10,10){$J_C$}}
\put(39,30){$g$}
\put(60,30){$f$}
\put(52,26){$q$}
\put(52,12){$p$}
\put(37,9){$g^*$}
\put(60,9){$f^*$}
\end{picture}
\end{center}
where $q:=(g,f)$,
$p=g^*+f^*$, $p\circ q=j_C$, and the inclusions
$F\hookrightarrow F\times B$, $B\hookrightarrow F\times B$
naturally identify $F$ with $F\times\{ 0\}$ and $B$ with
$\{0\} \times B$. Since $p:F\times B\rightarrow J_C$ 
is a finite covering of degree $r^2$, and since
$p^*(W_1)\equiv r\big( F\times\{ 0\} +\{0\} \times B\big)$
(cf \cite{Mumf2}; see also \cite {Ka2}), we 
get the relations:
$$
r^2\mbox{deg}(\widehat{c}\circ f^*)=\big( \{p^*(f^*B)\} . \{p^*S\} \big) ,\;
r^2\mbox{deg}(\widehat{c}\circ g^*)=\big( \{p^*(g^*F)\} . \{p^*S\} \big) ,$$
$$
r^2\mbox{deg}(\widehat{c}\circ j_C)=r\big( \{F\times\{ 0\} +
\{0\} \times B\} . \{p^*S\} \big) .
$$

Because 
$\{ p^*(g^*F+f^*B)\} .\{p^*S\}=
r^2\big(\{F\times \{ 0\} +\{0\}\times B\}.\{p^*S\}\big)$,
 the conclusion follows.

\begin{sectiune}
Proof of the Theorem.
\end{sectiune}

\noindent
In the sequel, we shall use the following simple
observation:

%
\begin{remark}
Let $E$ and $F$ be two elliptic curves and let
$f:F\rightarrow E$ be an isogeny. If $\mbox{deg}(f)$ and $r$
are coprime, then the morphism induced
between the $r$-torsion points, $f[r]:F[r]\rightarrow E[r]$
is an isomorphism.

Indeed, $\mbox{Ker}(f[r])\subset F[r]\cap \mbox{Ker}(f)$ and the
order of $F[r]$ and the order of $\mbox{Ker}(f)$
are both divisible by the order of $\mbox{Ker}(f[r])$.
It follows $\mbox{Ker}(f[r])=0$, and thus  $f[r]$ is an isomorphism.
\end{remark}

Let us denote $R:=r^2\Delta (r,c_1,c_2)=rc_2-(r-1)c_1^2/2$, 
and let $d$ be the greatest common divisor of the 
integers $r$ and $R$. 
The existence of a rank-$r$ holomorphic vector bundle
with Chern classes $c_1$ and $c_2$, in the particular case
$d=r$, was
proved in \cite{To2} by using unramified coverings
of $B$ with a suitable elliptic curve and  deformations
of sheaves.
 We can assume therefore that $d<r$, and we divide the proof of the
Theorem in several steps.\\\\
{\em Step 1.}
We reduce the proof to the case $d=1$. 
Suppose $d>1$. Since $R$ is divisible by $d$, it follows
that $c_1^2/2$ is divisible by $d$ as well. By means 
of the Lemma in \cite{To2}, there exists 
an unramified covering of degree $d$,
$t:B'\rightarrow B$ with a suitable elliptic curve
$B'$ such that, denoting by $\widetilde{t}:X':=
B'\times_BX\rightarrow X$ the canonically induced
unramified covering, where
$X'\stackrel{\pi '}{\rightarrow}B'$ is also a primary 
Kodaira surface, there exists
a class $c_1'\in NS(X')$ with $\widetilde{t}^*(c_1)=dc_1'$.
Set $r'=r/d$, $R'=R/d$, and $c_2'=c_2-(d-1){c_1'}^2/2$.
Then $R'=r'^2\Delta(r',c_1',c_2')$ and $R'$ and $r'$ are 
coprime.
If there exists a holomorphic rank-$r'$ vector bundle
${\cal E'}$ on the primary Kodaira surface $X'$ 
with Chern classes $c_1({\cal E'})=c_1'$ and $c_2({\cal E'})=c_2'$,
then we choose ${\cal E}:=\widetilde{t}_*({\cal E'})$.
A simple computation shows that ${\cal E}$ is a rank-$r$
holomorphic vector bundle on $X$ with Chern classes 
$c_1({\cal E})=c_1$ and $c_2({\cal E})=c_2$.\\\\
{\em Step 2.}
Suppose next $d=1$ and consider a cyclic isogeny
$\delta :F\rightarrow E$ of degree $R$, where $F$ is 
a suitable elliptic curve.
Let $\widehat{c}_1:B\rightarrow E$ be the isogeny
induced by the class $c_1\in NS(X)$. From Lemma
1.3 and Remark 3.1 we see that 
both 
morphisms $\delta [r]:F[r]\rightarrow E[r]$ and 
$\widehat{c}_1[r]:
B[r]\rightarrow E[r]$ are isomorphisms. Set $\psi:=
\delta [r]^{-1}\circ 
\widehat{c}_1[r]:B[r]\stackrel{\cong}{\longrightarrow}F[r]$.
From the properties of the Weil pairings (cf. \cite{Hu}, 12.2.4)
it follows that $\psi $ is an anti-isometry. Now, we
distinguish two cases:\\\\
{\em Case (a).} If $\psi$ is irreducible, denote by
$H_\psi=\mbox{Graph}(\psi)\subset (B\times F)[r]$.
Weil's Theorem (\cite{We} Satz 2; see also \cite{Ka2}, \cite{FrKa})
ensures that the quotient $J_\psi:=(B\times F)/H_\psi$ 
is the Jacobian $J_C$ of a curve $C$ of genus
2 which covers $r$-to-$1$ the elliptic curves $B$ and $F$.
Moreover, the morphism $\widetilde{c}:B\times F\rightarrow E$
defined by $\widetilde{c}(x,y)=\widehat{c}_1(x)-\delta(y)$
for any pair $(x,y)\in B\times F$ factors through
a morphism $\widehat{c}:J_C=(B\times F)/H_\psi\rightarrow E$.
It is clear that $\widehat{c}\circ f^*=\widehat{c}_1$ and
$\widehat{c}\circ g^*=\delta$.

We choose next a line bundle $L$ on the principal
elliptic bundle $Y:=X\times _BC\rightarrow C$, whose
Chern class modulo $\mbox{Tors}(NS(Y))$ equals $\widehat{c}$.
If $\varphi :Y\rightarrow X$ denotes the corresponding
covering of degree $r$, then Lemma 2.1 precisely says
that ${\cal E}:=\varphi _*L$ is a holomorphic rank-$r$
vector bundle with $c_1({\cal E})\equiv c_1 \mbox{modulo Tors}(NS(X))$
and discriminant 
$\Delta ({\cal E})=\mbox{deg}(\delta )/r^2=\Delta(r,c_1,c_2)$.\\\\
{\em Case (b).} If $\psi$ is reducible, we use
the "Reducibility criterion" from \cite{Ka2}.
Denote again $H_\psi=\mbox{Graph}(\psi)\subset (B\times F)[r]$
and take the quotient $J_\psi=(B\times F)/H_\psi$.
The morphism $\widetilde{c}:B\times F\rightarrow E$
defined by $\widetilde{c}(x,y)=\widehat{c}_1(x)-\delta(y)$
for any pair $(x,y)\in B\times F$ factors through
a morphism $\widehat{c}:J_\psi\rightarrow E$.
This time $J_\psi$ is no longer the Jacobian of a curve
of genus 2. Since $\psi$ is reducible, we get the
so-called ``diamond configuration'' (cf. \cite{Ka2}),
i.e. there exist two elliptic  curves $E_1$ and $E_2$,
an integer $k$ with $1\leq k<r$,
and isogenies 
$h:B\rightarrow F$, $f_i:B\rightarrow E_i$,
$f_i':E_i\rightarrow F$, for $i=1,2$, such that
$h=f_1'\circ f_1=f_2'\circ f_2$, $\mbox{deg}(f_1)=\mbox{deg}
(f_2')=r-k$, $\mbox{deg}(f_1')=\mbox{deg}(f_2)=k$.
Moreover, if $\widetilde{f_i'}:F\rightarrow E_i$ denotes
the dual map of $f_i'$, for $i=1,2$, then
$J_\psi$ is isomorphic to  $E_1\times E_2$
via the map $p:J_\psi\rightarrow E_1\times E_2$,
$p(\widehat{x,y}):=\big( f_1(x)-\widetilde{f_1'}(y),
f_2(x)+\widetilde{f_2'}(y) \big)$.

The natural inclusion of $B$ in $E_1\times E_2$ is given by
$(f_1,f_2)$ and the natural inclusion of $F$ in $E_1\times E_2$
is given by $(-\widetilde{f_1'},\widetilde{f_2'})$.
It follows easily that $\widehat{c}\circ (f_1,f_2)=
\widehat{c}_1$
and $\widehat{c}\circ (-\widetilde{f_1'},\widetilde{f_2'})=-\delta$.
 We consider the maps
$c':E_1\rightarrow E$
and $c'':E_2\rightarrow E$ given by
$c'(u)=\widetilde{c}(u,0)$, and $c''(v)=\widetilde{c}(0,v)$,
respectively. Then a simple computation
gives $\widetilde{c}_1=c'\circ f_1+c''\circ f_2$.
Moreover, composing the 
relation $\widehat{c}\circ (-\widetilde{f_1'},\widetilde{f_2'})=-\delta$
by $h$ to the right we get:
\begin{equation}
k(c'\circ f_1)-(r-k)(c''\circ f_2)=\delta \circ h.
\end{equation}
Consider now the unramified coverings
$\widetilde{f_1}:E_1\rightarrow B$ and
$\widetilde{f_2}:E_2\rightarrow B$
of degree $(r-k)$ and $k$ respectively.
Let $X_1\rightarrow E_1$ and $X_2\rightarrow E_2$
be the primary Kodaira surfaces defined by these
coverings by taking the fibered products. Choose line
bundles $L_1$ and $L_2$ on $X_1$ respectively on $X_2$,
with Chern classes modulo torsion $c'$ and $c''$ respectively.
Denote $X_1\stackrel{\varphi _1}{\rightarrow }X$
and $X_2\stackrel{\varphi _2}{\rightarrow }X$ 
the induced unramified coverings of degrees
$(r-k)$, respectively $k$, and take the holomorphic
vector bundles on $X$, ${\cal E}_1:=\varphi_{1,*}(L_1)$,
and ${\cal E}_2:=\varphi_{2,*}(L_2)$, of ranks 
$(r-k)$, respectively $k$. These vector bundles
have Chern classes (modulo torsion) $c'\circ f_1$, respectively,
$c''\circ f_2$, and vanishing discriminants.
Finally, set ${\cal E}:={\cal E}_1\oplus {\cal E}_2$.
Then ${\cal E}$ is a rank-$r$ vector bundle
with $c_1({\cal E})\equiv c_1$ modulo $\mbox{Tors}(NS(X))$.
Now, using (3.1) a simple computation
leads us to :
$$\Delta ({\cal E}_1\oplus {\cal E}_2)=\frac{1}{2r}\left[
\frac{\widehat{c}_1^2}{r}-\frac{(c'\circ f_1)^2}{r-k}
-\frac{(c''\circ f_2)^2}{k}\right] =\Delta (r,c_1,c_2).$$
\\\\
{\em Step 3.}
In order to get rid of the torsion,
we need to add multiples of a class of a fiber.
To do this, in the irreducible case,
we consider $F_Y$ to be a fiber
of the projection map $Y\longrightarrow C$,
$F_X$ to be a fiber
of the projection map $X\longrightarrow B$,
and for any $m\geq 0$, we set $L_m={\cal O}_Y(mF_Y)$.
In homology, $\varphi _*\{ F_Y\} =\{ F_X\} $,
and thus $\varphi _*\big( c_1\left(
L_m\right) \big) =mc_1\left( {\cal O}_X(F_X)\right) $. Formula
(2.2) reads here $c_1(\varphi _*(L\otimes L_m))
=c_1({\cal E})+mc_1\left( {\cal O}_X(F_X)\right) $.
In the reducible case, we apply a similar argument,
using a fiber of $X_1\rightarrow X$ or of $X_2\rightarrow X$, and this
ends the proof.

\begin{remark}
For a primary Kodaira surface the compactness
theorem in \cite{To4} (Theorem 5.9) combined with
the existence result above, produces moduli
spaces of stable bundles, which are non-empty, holomorphically
symplectic compact manifolds, when the Chern classes are chosen
in the 
stably irreducible range as in \cite{To4}. 
For example, if $c_1\in NS(X)$ is chosen such that
$m(2,c_1)>0$ and 
$0\le \Delta (2,c_1,c_2)<m(2,c_1)$,
we are in 
this
range and all the 2-vector bundles
with 
the 
given invariants are
%
stable
with respect to any Gauduchon metric on the primary Kodaira surface.
\end{remark}


\begin{thebibliography}{abcdefg}

\bibitem[BaL]{BaL} B\u anic\u a C.,  Le Potier, J.: Sur l'existence des
fibr\'es
%
vectoriels holomorphes sur les surfaces non-alg\'ebriques, J. Reine Angew.
Math.
%
\textbf{378}, (1987) 1-31
\bibitem[BPV]{BarPV}  Barth, W., Peters, C., Van de Ven, A.: Compact
complex
surfaces, Springer-Verlag: Berlin-Heidelberg-New York, 1984
\bibitem[Bri1]{Bri1} Br{\^{\i}}nz{\u{a}}nescu, V.: N\' eron-Severi group for
non-algebraic elliptic surfaces I: elliptic bundle case, Manuscripta Math.
{\bf 79}, (1993) 187-195; II: non-K\" ahlerian case,  Manuscripta Math.  
{\bf 84}, (1994) 415-420; III, Rev. Roumaine Math. Pures
Appl. {\bf 43} (1998), 1-2, 89-95
\bibitem[Bri2]{Bri2} Br\^{\i}nz\u{a}nescu, V.: The Picard group of a 
primary Kodaira surface, Math. Ann., {\bf 296}, (1993) 725--738
\bibitem[Bri3]{Bri3}  Br\^{\i}nz\u{a}nescu, V.:
Holomorphic vector bundle over compact
complex surfaces,  Lect. Notes in Math. {\bf 1624}, Springer 1996
\bibitem[BrF]{BrF} Br{\^{\i}}nz{\u{a}}nescu, V., Flondor, P.:
Holomorphic
2-vector bundles on non-algebraic 2-tori, J. reine angew. Math. {\bf 363},
(1985) 47--58
\bibitem[ElFo]{ElFo} Elencwajg, G., Forster, O.: Vector bundles on
manifolds
without divisors and a theorem of deformation, Ann. Inst. Fourier {\bf
32}(4), (1982) 25--51
\bibitem[FrKa]{FrKa}  Frey, G., Kani, E.: Curves of genus 2
covering elliptic curves and an arithmetical application,
Arithmetic Algebraic Geometry, Progress in Math {\bf 89} (1991),
153-176
\bibitem[Hf]{Hf} H{\"{o}}fer, T.: Remarks on torus principal bundles,
J. Math. Kyoto Univ. (JMKYAZ) {\bf 33}(1), (1993) 227--259
\bibitem[Hu]{Hu} Husem\" oller, D.: Elliptic curves, Graduate Texts in
Math. {\bf 111}, Springer-Verlag 1987
\bibitem[Ka1]{Ka1}  Kani, E.:  Elliptic curves on abelian
surfaces, Manuscripta Math. {\bf 84} (1994), 199-223
\bibitem[Ka2]{Ka2}  Kani, E.: The number of curves of genus two with
elliptic differentials, J. reine angew. Math. {\bf 485},
(1997) 93-121                     
\bibitem[Kod]{Kod}  Kodaira, K.: On the structure of compact complex
analytic
surfaces I , Amer. J. Math. {\bf 86}, (1964) 751--798
\bibitem[Kr]{Kr}  Krazer, A.: Lehrbuch der Thetafunktionen, Leipzig 1903 
(Chelsea Reprint, 1970)
\bibitem[LanBi]{LanBi} Lange, H., Birkenhake, Ch.:
Complex Abelian Varieties.
Grund. Math. Wissenschaften,
Springer-Verlag (1992)
\bibitem[LeP]{LeP} Le Potier, J.: Fibr\'{e}s vectoriels sur les surfaces 
K3, S\'{e}minaire Lelong-Dolbeault-Skoda, LNM {\bf 1028}, Springer 1983
\bibitem[Mumf1]{Mumf1}  Mumford, D.: Abelian varieties, Oxford
University Press, 1974
\bibitem[Mumf2]{Mumf2}  Mumford, D.: Prym Varieties I,
Contribution to Analysis, Acad. Press, New York (1974) 325-350
\bibitem[Schw]{Schw} Schwarzenberger, R.L.E.: Vector bundles on algebraic
surfaces,
Proc. London Math. Soc. {\bf 3}, (1961) 601-622
\bibitem[Te]{Te} Teleman, A.: Moduli Spaces of Stable Bundles on non-K\"ahlerian
Elliptic Fibre Bundles over Curves, Expo. Math., {\bf 16}, (1998), 193-248
\bibitem[To1]{To1} Toma, M.: Une classe de fibr\'es vectoriels
holomorphes
sur
les 2-tores complexes, C.R. Acad. Sci. Paris, {\bf 311}, (1990)
257-258
\bibitem[To2]{To2} Toma, M.: Stable bundles on non-algebraic surfaces
giving
rise to compact moduli spaces, C.R. Acad. Sci. Paris, {\bf 323},
(1996) 501-505
\bibitem[To3]{To4} Toma, M.: Stable bundle with small
second Chern classes on surfaces, Heft 209 (1999), Preprint
Univ. Osnabr\"uck
\bibitem[We]{We} Weil, A.: Zum Beweis des Torellischen Satzes, 
G\"{o}ttinger Nachr., no. 2, (1957) 33--53

\end{thebibliography}
\end{document}